\documentclass[11pt,leqno]{paper}
\usepackage{amsmath,amssymb,amsfonts}
\usepackage{graphics}
\usepackage[dvips]{epsfig}
\usepackage{subfigure}
\usepackage{color}
\usepackage{hyperref}
\hypersetup{
    colorlinks,
    citecolor=black,
    filecolor=black,
    linkcolor=black,
    urlcolor=black
}

 \newtheorem{theorem}{{Theorem}}[section]
\newtheorem{isom.ext}[theorem]{{Trivial isometric extension}}%[section]
\newtheorem{lemma}[theorem]{{Lemma}}%[section]
\newtheorem{corollary}[theorem]{{Corollary}}%[section]
\newtheorem{remark}[theorem]{{Remark}}%[section]
\makeindex

\title{ Torsional Rigidity on Compact Riemannian Manifolds with lower Ricci Curvature Bounds}

\author{Najoua Gamara,\; Abdelhalim Hasnaoui, \;Akrem Makni  \\
{\small\em  \today }}

 \date{ }
\begin{document}
\maketitle
%\dominitoc
\makeatletter

\def\thickhrulefill{\leavevmode \leaders \hrule height 1ex \hfill \kern \z@}

\vspace*{1cm} {\bf \textbf{Abstract}:}

In this article we prove a reverse H\"older inequality for the fundamental eigenfunction
 of the Dirichlet problem on domains of a compact Riemannian manifold with lower
  Ricci curvature bounds. We also prove an isoperimetric inequality for the torsional ridigity of such domains.\\

$Mathematics \; Subject \; Classification \; (2000)$: 35P15, 60J65, 53C21, 58J60.\\
$Key \; words$: isoperimetric inequalities, eigenfunctions, Ricci
curvature, Reverse H\"older Inequality, torsional rigidity.
 %\tableofcontents

\section{Introduction and Main Results}\label{ricci}
%\subsection{Motivation}

In $1972,$  following the spirit of the works of Faber and Krahn \cite{krahn,faber},  Payne and Rayner \cite{PR1,PR2}
    proved a reverse H\"older inequality for the norms $L_{1}$ and $L_{2}$ of the first
   eigenfunction of the Dirichlet problem on bounded  domains $D$ of
   $\mathbb{R}^{2}:$
$$\Big(\int_{D}\,u_{1}\,dxdy\Big)^{2}\geq\,\displaystyle\frac{4\pi}{\lambda_{1}}\int_{D}\,u_{1}^{2}\,dxdy$$
where $\lambda_{1}$ is the lowest eigenvalue of the fixed membrane
problem, and $u_{1}$  the corresponding eigenfunction. Equality holds if and only if $D$ is a disc.\\
   The work of Payne and Rayner was generalized by K\"ohler and Jobin
    \cite{KJ1} for bounded domains of  $\mathbb{R}^{n},$  $n\geq3.$
   In 1982, G.Chiti \cite{Ch1}, generalized the reverse H\"older inequality for the norms $L_{q}$
    and $L_{p},$  $q\geq p >0$ for bounded domains of $\mathbb{R}^{n},$ $n\geq2.$ \\
$$(\int_{D}u^{q}dxdy)^{\frac{1}{q}} \leq\,K(p,q,\lambda,n)\,(\int_{D}u^{p}dxdy)^{\frac{1}{p}}$$
Where
$$K(p,q,\lambda,n)=(nC_{n})^{\frac{1}{q}-\frac{1}{p}}\lambda^{\frac{(q-p)n}{2pq}}
j_{\frac{n}{2}-1,1}^{n(\frac{1}{q}-\frac{1}{p})}\displaystyle\frac{\Big(\int_{0}^{1}r^{n-1+q(1-\frac{n}{2})}J_{\frac{n}{2}-1}^{q}(j_{\frac{n}{2}-1,1}\,r)dr\Big)^{\frac{1}{q}}
}{\Big(\int_{0}^{1}r^{n-1+p(1-\frac{n}{2})}J_{\frac{n}{2}-1}^{p}(j_{\frac{n}{2}-1,1}\,r)dr\Big)^{\frac{1}{p}}}$$
This inequality is isoperimetric in the sens that equality holds if and only if $D$ is a ball.

      The main ideas of this paper  were  early investigated in the PHD thesis of H.Hasnaoui \cite{H.Has}, the first one was  to establish a generalization of  Chiti's reverse H\"older inequality for the norms $L_{q}$
       and $L_{p},$  $q\geq p >0$ for Compact Riemannian manifolds with Ricci curvature bounded from below and the second one is a version of the Saint Venant Theorem for such manifolds.
      \par   In fact Modern geometric analysts, including
Chavel and Gromov, have identified such manifolds as important, and have
related the Ricci bound to many estimates of eigenvalues, as well as to other quantities
of interest in differential geometry. Thus it is very natural to consider similar questions
about the Laplacian and the torsional rigidity in this context. Both are isoperimetric
results, is that in both cases the quantity of interest is dominated by the analogous
expression on spheres.
 Much of the background and many results for the spectrum of such manifolds could be found in \cite{Cheng}, \cite{BM}, \cite{Gam} and more recently in \cite{LingLu}.
In 1856, Saint-Venant \cite{Saint-Venant} observed that columns with circular
cross-sections offer the greatest resistance to torsion for a given
cross-sectional area. This fact was proved a century or so later by
P\'olya using Steiner symmetrization \cite{polya48}. We also note
the independent proof  by Makai in 1963 one can see \cite{Makai}.
Torsional rigidity is a physical quantity of much interest, see for
example \cite{Bandle, PolyaSzego}, the  recent
\cite{RatzkinCaroll2011,Iversen}, the more recent \cite{vdb0,RatzkinCaroll2014} and the classical papers of
Payne \cite{Payne1,Payne0}, Payne-Weinberger
\cite{PayneWeinbergerTorsion}. As we show here, it is also a
quantity of much interest for geometers as
well.\\
Parts of our work  follow an analysis of Ashbaugh and Benguria
\cite{AshbaughBenguriaSphere} for subdomains of hemispheres, one can see also\cite{ChavelFeldman}.
In the sequel we will introduce the main results of this paper:\\

Let $(M, g)$ be a compact connected Riemannian
manifold of dimension $n \geq 1$ without boundary. We denote by
$$R(M,g) = \inf \left\{ Ric(v,v), \;\;\;\; v\in UT(M)  \right\},$$
the infimum of the Ricci curvature $Ric$ of $(M, g),$ here $UT(M)$
is the unit tangent fiber bundle of the manifold $(M, g)$. Let
$(\mathbb{S}^{n},g^{\star})$ be unit the sphere of the space
$\mathbb{R}^{n+1}$, endowed with the induced metric, then
$R(\mathbb{S}^{n},g^{\star}) = n-1$. We suppose, as done in
\cite{BM}, that $R(M,g)$ is strictly positive, and normalize the
metric $g$ so that $R(M, g)\geq R(\mathbb{S}^{n},g^{\star}) = n-1$.
In the sequel, we will denote by $V(M)=\int dv_{g}$ the volume of
$(M, g),$ where the element of volume is denoted $dv_{g}$,  $\omega_{n}$ the volume of the unit
sphere
$(\mathbb{S}^{n},g^{\star})$ and  $\beta=V(M)/\omega_{n}.$\\
Let $D$ be a connected bounded domain of $M$  with
smooth boundary, and $D^{\star}$ the geodesic ball of
$\mathbb{S}^{n}$ centered at the north pole such that
$$Vol (D)=\beta \, Vol(D^{\star}).$$

  We are interested in the comparison of the fundamental solutions
  $u$ and $v$
  of problems $(P_{1})$ and $(P_{2})$ respectively:

\begin{eqnarray*}(P_{1})
\left\{\begin{array}{llll}
\Delta u + \lambda u &=& 0 &\text{ in }  D\\
u &=&0&   \text{ on } \partial D,
\end{array}
\right.
\end{eqnarray*}

\begin{eqnarray*}(P_{2})
\left\{\begin{array}{llll}
\Delta v + \lambda v &=& 0 &\text{ in }  B_{\theta_{1}(\lambda)}\\
v &=&0&   \text{ on } \partial B_{\theta_{1}(\lambda)},
\end{array}
\right.
\end{eqnarray*}
 $B_{\theta_{1}(\lambda)}$ is  the geodesic ball of $\mathbb{S}^{n}$
centered at the north pole of radius $\theta_{1}=\theta_{1}(\lambda),$ where $\lambda$ is the first
eigenvalue for the Dirichlet  problem $(P_{2}),$ and $\Delta$
denotes indifferently the laplacian operator on $M$ or
$\mathbb{S}^{n}.$\\
Let $u^{\ast}$ be the decreasing rearrangement of $u$ and $u^{\star}$ the
corresponding radial function defined on $D^{\star}$ the geodesic
ball of $(\mathbb{S}^{n},g^{\star})$ which has the same relative
volume as $D,$ (see section2 for notations and details).

We prove the following result
\begin{theorem}
 Let $u$ be the fundamental eigenfunction of problem $(P_{1})$. Let $p> 0$ and $v$ be the fundamental eigenfunction of  problem $(P_{2})$ chosen such that
 \begin{eqnarray}
\int_{D} u^{p} dv_{g} =\beta \, \int_{B_{\theta_{1}(\lambda)}} v^{p}
dv_{g^{\star}}.
\end{eqnarray}
then, there exists $\theta_{2} \in (0,\theta_{1}),$ such that
\begin{eqnarray}
% \nonumber to remove numbering (before each equation)
 u^{\star}(\theta) \leq v( \theta) , & \hbox{$\forall \theta \in [0, \theta_{2}] $ ;} \\
    u^{\star}(\theta) \geq v(\theta), & \hbox{$\forall  \theta \in [\theta_{2}, \theta_{1}] $.}
\end{eqnarray}
 \end{theorem}
As a consequence of Theorem 1.1, we obtain:
\begin{corollary}: Chiti's Reverse H\"older Inequality for Compact Manifolds\\
  Let  $p , q$ be real numbers such that $q \geq p > 0$, then $u$ and $v$ are related by this inequality
\begin{eqnarray}
\frac{(\int_{D} u^{q} dv_{g} )^{\frac{1}{q}}}{(\int_{D} u^{p} dv_{g}
)^{\frac{1}{p}}} \leq \beta^{\frac{1}{q}-\frac{1}{p}} \,
\frac{(\int_{B_{\theta_{1}(\lambda)}}v^{q}
dv_{g^{\star}})^{\frac{1}{q}}}{(\int_{B_{\theta_{1}(\lambda)}}v^{p}dv_{g^{\star}})^{\frac{1}{p}}},
\end{eqnarray}
with equality if and only if the triplet $(M,D, g)$ is isometric to
the triplet $(\mathbb{S}^{n}, D^{\star},g^{\star})$.
\end{corollary}
Next, we focus on  the torsional rigidity $T(D)$ of the manifold
$D.$ Recall that:
 \begin{eqnarray}
T(D) = \int_{D} w dv_{g},
\end{eqnarray}
 where $w$ is is the  smooth solution of the boundary value problem of Dirichlet-Poisson type ( $w$ is called the warping function of $D$)
  \begin{eqnarray} \label{ss}
\Delta_{g} w + 1 & = & 0 \mbox{ in } D, \nonumber \\
w & = & 0 \mbox{ on }\partial D.
\end{eqnarray}
We obtain the following result:

\begin{theorem}: Saint Venant Theorem for Compact Manifolds\\
 Let $(M, g)$ be a complete  Riemannian manifold of dimension $n$,
without boundary satisfying $R(M, g)\geq  n-1$, and $D$ a sub
manifold of $M$ of the same dimension with smooth boundary, we have
 \begin{eqnarray}
T(D)  \leq \beta \, T(D^{\star}).
\end{eqnarray}
 the equality holds if and only
 if the triplet $(M,D,g)$ is isometric to the triplet $(\mathbb{S}^{n},D^{\star},g^{\star})$.
\end{theorem}

And finally we give a comparison formula for the warping function
$w$ which allows us to obtain directly the result of theorem1.3.

\section{Preliminary Tools}
Denote by $\lambda (D)$ the first eigenvalue of the Laplacian
$\Delta$ for the Dirichlet problem  on $D$, and let $u$ be the
positive associated eigenfunction. Therefore $u$ satisfies
\begin{eqnarray*}(P_{1})
\left\{\begin{array}{llll}
\Delta u + \lambda u &=& 0 &\text{ in }  D\\
u &=&0&   \text{ on } \partial D,
\end{array}
\right.
\end{eqnarray*}
The variational formula of $\lambda (D)$ is given by
  $$\lambda (D) = \inf \left\{ \frac{\int_{D}|df|^{2} \, dv_{g}}{\int_{D}|f|^{2} \, dv_{g}};\;\;\;f \ncong 0; \;\;\;C^{1}- piecewise \;\; ;\;\;\; f_{|_{\partial D}} \equiv 0    \right\},$$
where $|df|$ is the Riemannian norm of the differential of $f$. We
have equality if and only if $f$ is of class $C^{2}$ and is an
eigenfunction associated with the first eigenvalue $\lambda (D).$

 For $0 \leq  t \leq \overline{u} = \sup u$, let $D_{t} =\big\{x \in D \big| u(x)>t\big\}$. Define the function
\begin{eqnarray}
V(t) = \int_{D_{t}} dv_{g}
\end{eqnarray}
The co-area formula gives
\begin{eqnarray}
V(t) =\int_{D_{t}} dv_{g} = \int_{t}^{\overline{u}}\int_{\partial
D_{\tau}}\frac{1}{|\nabla u|} d\sigma \; d\tau.
\end{eqnarray}
Here $d\sigma$ is the (n-1)-dimensional Riemannian measure in
$(M,g)$. For what follows, we will also denote the
$(n-1)$-dimensional Riemannian measure in
$(\mathbb{S}^{n},g^{\star})$ by $d\sigma$. Since $D$ has bounded
measure, the above shows that the function
\begin{eqnarray}
  t \mapsto \int_{\partial D_{t}}\frac{1}{|\nabla v|} d\sigma
\end{eqnarray}
is integrable, and therefore the function $V$ is absolutely
continuous. Hence $V$ is differentiable almost everywhere and
\begin{eqnarray}
\frac{d V}{d t} = - \int_{\partial D_{t}}\frac{1}{|\nabla u|}
d\sigma \; < \; 0
\end{eqnarray}
for almost all $t \in [ 0, \overline{u}]$. The function $V$ is then
a  non increasing function and has an inverse which we denote by
$u^{\ast}.$ The function $u^{\ast}$ is absolutely continuous. Now,
applying the Cauchy-Schwartz inequality, we obtain
\begin{eqnarray}
  \left(\int_{\partial D_{t}} d\sigma\right)^{2} \leq \left(\int_{\partial D_{t}}\frac{1}{|\nabla u|} d\sigma\right)\, \left(\int_{\partial D_{t}}|\nabla u|d\sigma\right).
\end{eqnarray}
Hence
\begin{eqnarray}
 -\frac{d u^{\ast}(s)}{ds} = - \frac{1}{V'(u^{\ast}(s))} \leq \frac{\int_{\partial D_{u^{\ast}(s)}}|\nabla u|d\sigma}{ \left(\int_{\partial D_{u^{\ast}(s)}}d\sigma\right)^{2}}.
\end{eqnarray}

We now use the following isoperimetric inequality due to M. Gromov
\cite{Gromov} which relates the volume of the boundaries of $D$ and
$D^{\ast}.$
\begin{lemma}  \label{lem3.5}
Under the same hypothesis given above
\begin{equation}
Vol_{n-1}(\partial D) \geq \beta \, Vol_{n-1}(\partial D^{\ast}),
\end{equation}
where  $Vol_{n-1}$ is the $(n-1)$-dimensional volume relative to $g$
and
$g^{\ast}$. Equality holds if and only if the triplet  $(M, D, g)$ is isometric to
the triplet $(\mathbb{S}^{n}, D^{\star}, g^{\star}).$
\end{lemma}
Let
\begin{eqnarray}
A(\theta) \equiv s =    \beta \, \omega_{n-1} \int_{0}^{\theta}
(\sin \tau)^{n-1} d\tau
\end{eqnarray}
The quantity $\omega_{n-1} \, \int_{0}^{\theta} (\sin \tau)^{n-1}
d\tau$ is the $n$-volume of the geodesic ball of radius $\theta$ in
$\mathbb{S}^{n}$.  If we let the function $L(\theta)$ denotes the
$(n-1)$-dimensional volume of the geodesic ball of radius $\theta$,
i.e.,
\begin{eqnarray}
L(\theta) = \omega_{n-1}\, (\sin\theta)^{n-1} =
\frac{A'(\theta)}{\beta}.
\end{eqnarray}
Inequality (14) can then be written as
\begin{eqnarray}
Vol_{n-1}(\partial D) \geq \beta \, L(\theta (Vol(D ))) =\beta \,
\omega_{n-1}(\sin\theta(Vol(D )))^{n-1}
\end{eqnarray}
where $\theta(s)$ is the inverse function of $A$ defined in (15).\\
Now, applying  inequality (17) to the domain $D_{u^{\ast}(s)}$ and
combining it with inequality (13), we obtain
\begin{eqnarray}
-\frac{d u^{\ast}(s)}{ds} \leq  (\beta \, \omega_{n-1})^{-2}(\sin
\theta (s))^{2-2n} \int_{\partial D_{u^{\ast}(s)}}|\nabla u|d\sigma.
\end{eqnarray}
We then apply Gauss Theorem to the Dirichlet problem on $D_{t}$,
using the smoothness of its boundary  $\partial D_{t},$ we get
\begin{eqnarray}
\int_{\partial D_{t}}|\nabla u| \; d\sigma = \int_{D_{t}}\Delta u \;
dv_{g}= \lambda \int_{D_{t}}u \; dv_{g},
\end{eqnarray}
here we used the fact that the outward normal to $D_{t}$ is given by
$-\frac{\nabla u}{|\nabla u|}$.\\
\textbf{ Remark} $\forall p \geq
0,$ we have
\begin{eqnarray}
 \int_{D_{u^{\ast}(s)}}u^{p} \; dv_{g} = \int_{u^{\ast}(s)}^{\overline{u}} \tau^{p}\int_{\partial D_{\tau}}\frac{1}{|\nabla u|} \; d\sigma d \tau = -\int_{u^{\ast}(s)}^{\overline{u}} \tau^{p} V'(\tau) d\tau.
\end{eqnarray}
The change of variables $\eta = V(\tau)$  gives
\begin{eqnarray}
\int_{D_{u^{\ast}(s)}}u^{^p} \; dv_{g}=
\int_{0}^{s}(u^{\ast}(\eta))^{p} d\eta.
\end{eqnarray}

Finally combining (20) for $p = 1$ with
equalities (18) and (19), we obtain the following.
\begin{lemma}
Let $u$ be a solution of problem $(P_{1})$. Then $u^{\ast}$, its
decreasing rearrangement, satisfies the integro-differential
  inequality
\begin{eqnarray}
 -u^{\ast \prime}(s) \leq  \lambda (\beta \, \omega_{n-1})^{-2}(\sin \theta(s))^{2-2n}\int_{0}^{s} u^{\ast}(\xi)d\xi.
\end{eqnarray}
\end{lemma}
for almost every $s > 0$.\\
Let  $B_{\theta_{1}(\lambda)}$ be the geodesic ball of $\mathbb{S}^{n}$ centered
at the north pole with radius $\theta_{1}=\theta_{1}(\lambda),$ such that the following problem
\begin{eqnarray*}(P_{2})
\left\{\begin{array}{llll}
\Delta v + \lambda v &=& 0 &\text{ in }   B_{\theta_{1}(\lambda)}\\
v &=&0&   \text{ on } \partial  B_{\theta_{1}(\lambda)}
\end{array}
\right.
\end{eqnarray*}
has a solution. Let
  $v > 0$ be the first Dirichlet
eigenfunction on $B_{\theta_{1}(\lambda)}$  then by lemmas 3.1 and 3.2 in
\cite{Benguria2011}, we conclude that $v$ depends only of $\theta$ and
is strictly decreasing on $[0,\theta_{1})$. Therefore we denote by
$v(\theta)$ the function $v$. So in polar coordinates the problem
$(P_{2})$ can be rewritten as
\begin{eqnarray}
 -((\sin \theta)^{n-1}v'(\theta))' = \lambda (\sin \theta)^{n-1} v(\theta)
\end{eqnarray}
in $[0, \theta_{1}]$. The boundary conditions are $v(0)$ finite and
$v(\theta_{1}) = 0$. Let $v^{\ast}$ be the function defined in $[0,
A(\theta_{1})]$ from the relation
\begin{eqnarray}
 v^{\ast}(A(\theta)) = v(\theta); \;\;\;\;\;\;  \forall \theta \in [0, \theta_{1}]
\end{eqnarray}
Integrating the equality (23), we obtain
\begin{eqnarray}
-(\sin \theta)^{n-1}v'(\theta) = \lambda \int_{0}^{\theta}(\sin
\alpha)^{n-1} v(\alpha)d\alpha
\end{eqnarray}
 Then using (24),  we can rewrite  the left-hand side of (25) as
\begin{eqnarray}
-(\sin \theta(s))^{n-1}v'(\theta(s)) = - \beta \, \omega_{n-1}(\sin
\theta(s))^{2n-2}v^{\ast \prime }(s)
\end{eqnarray}
for all $s$ in $[0, A(\theta_{1})]$. The change of variables $\zeta
= A(\alpha)$ in  the right-hand side of (25) gives
\begin{eqnarray}
\lambda \int_{0}^{\theta(s)}(\sin \alpha)^{n-1} v(\alpha)d\alpha =
\lambda (\beta \, \omega _{n-1})^{-1}\int_{0}^{s}
v^{\ast}(\zeta)d\zeta
\end{eqnarray}
Finally from (26) and (27),  we obtain
\begin{eqnarray}
-v^{\ast \prime}(s) = \lambda (\beta \, \omega _{n-1})^{-2}(\sin
\theta(s))^{2-2n}\int_{0}^{s} v^{\ast}(\xi)d\xi
\end{eqnarray}

\section{Chiti's Reverse H\"older Inequality and the Saint-Venant Theorem for
Compact Manifolds}

\subsection{Chiti's Reverse H\"older Inequality}

 In this section we will
 prove the extension of the Chiti's Comparison Lemma, given for
domains in the case of $\mathbb{R}^2$ and $\mathbb{R}^{n}$ in the
original papers of Payne-Rayner \cite{PR1,PR2}, then extended by
Kohler-Jobin \cite{KJ1} and Chiti \cite{Ch1, Ch2}. In
\cite{AshbaughBenguriaSphere} a very general version of the Chiti's
Comparison Lemma is available for domains of the sphere
$\mathbb{S}^{n-1}$. Our work deals with the case of a complete
Riemannian manifold with bounded Ricci
curvature.\\
 Let $u^{\star}$  be the radial function defined in $D^{\star}$ by
\begin{eqnarray}
u^{\star}(\theta) = u^{\ast}(A(\theta))
\end{eqnarray}
\textbf{CLAIM}: we have
\begin{eqnarray} \label{clm}
vol (D)\geq \beta \, Vol(B_{\theta_{1}(\lambda)} )
\end{eqnarray}
\textbf{Proof of the Claim}:  Assume in the contrary that $\beta \, Vol(B_{\theta_{1}(\lambda)} ) > vol (D),$  since
 $vol(D) = \beta\, Vol(D^{\star}),$ we obtain
$Vol(B_{\theta_{1}(\lambda)} )> vol(D^{\star}).$ Then using the fact that the two
geodesic balls are centered at the north pole of $\mathbb{S}^{n},$
we obtain $D^{\star} \subsetneq B_{\theta_{1}(\lambda)} $. Finally from  domains
monotonicity of eigenvalues, we deduce  that $\lambda_{1}(D^{\star})
> \lambda_{1}(B_{\theta_{1}(\lambda)}) = \lambda$ which contradicts the result
of  Theorem 5 in \cite{BM}.\\ Now, we will  prove the following result which is crucial for
the proof of theorem1.1

\begin{lemma}  \label{lem3.5}
 Let $v$ be the solution of $(P_{2})$ chosen such that $v(0)= u^{\star}(0)$. Then
\begin{eqnarray}
v(\theta) \leq u^{\star}(\theta), \;\;\;\;\; \forall \theta \in [0,
\theta_{1}] .
\end{eqnarray}
\end{lemma}
\textbf{Proof} Using  (30),  two cases occur:\\
\textbf{First Case:}
$vol(D) = \beta \, Vol(B_{\theta_{1}(\lambda)}).$ Then we have $vol(B_{\theta_{1}(\lambda)})
= vol(D^{\star})$, hence
\begin{eqnarray}
\int_{B_{\theta_{1}(\lambda)}} |\nabla u^{\star}|^{2} dv_{g^{\star}} =
\omega_{n-1}\, \int_{0}^{\theta_{1}} \left(\frac{d
u^{\star}(\theta)}{d \theta}\right)^{2} \sin^{n-1}\theta d \theta
\end{eqnarray}
Using the change of variables $A(\theta) = s$, we obtain
\begin{eqnarray}
 \int_{0}^{\theta_{1}} \left(\frac{d  u^{\star}(\theta)}{d \theta}\right)^{2} \sin^{n-1}\theta d \theta = \beta \, \omega_{n-1} \int_{0}^{A(\theta_{1})} (u^{\ast '}(s))^{2} (\sin\theta(s))^{2n -2}  ds
\end{eqnarray}
Then we combine equalities (32) and (33) with lemma2.2 to get
\begin{eqnarray}
 \int_{B_{\lambda}} |\nabla u^{\star}|^{2} dv_{g^{\star}} \leq  \lambda \, \beta^{-1} \, \int_{0}^{A(\theta_{1})} ((-u^{\ast '}(s)\int_{0}^{s}u^{\ast}(\xi) d\xi)ds
\end{eqnarray}
Now, an integration by parts in the second member of the last
inequality and the fact that
 \begin{eqnarray}
 \int_{B_{\theta_{1}(\lambda)}}  u^{\star 2} dv_{g^{\star}} = \beta^{-1} \, \int_{0}^{A(\theta_{1})}(u^{\ast}(s))^{2} ds
\end{eqnarray}
gives
\begin{eqnarray}
\frac{\int_{B_{\theta_{1}(\lambda)}} |\nabla u^{\star}|^{2}
dv_{g^{\star}}}{\int_{B_{\theta_{1}(\lambda)}}  u^{\star 2} dv_{g^{\star}}} \leq
\lambda
\end{eqnarray}
Using that $\lambda$ is the minimum of the Rayleigh quotient on $B_{\theta_{1}(\lambda)}$, it follows that this minimum is achieved for $u^{\star}$,
 and so $u^{\star}$ is indeed an eigenfunction associated with $\lambda$ on $B_{\theta_{1}(\lambda)}$. Now, using the simplicity of the fundamental
eigenvalue and the hypothesis of our Lemma, we get $u^{\star} = v$.\\
\textbf{Second Case:}  $vol(D) > \beta \, Vol(B_{\theta_{1}(\lambda)})$.\\
 In one hand, we have
$u^{\ast}(A(\theta_{1})) > 0 $ while $v^{\ast}(A(\theta_{1})) = 0$.
In the other hand
\begin{eqnarray}
     u^{\ast}(0)=u^{\star}(0) = v(0) = v^{\ast}(0)
\end{eqnarray}
 then, we can find a constant $c \geq 1,$ such that
\begin{eqnarray}
c \,  u^{\ast}(s)   \geq v^{\ast}(s) \quad \forall s \in
[0,A(\theta_{1})]
\end{eqnarray}
Let $c'$ be the constant defined by
\begin{eqnarray}
c' = \inf\{ c \geq 1\; ; \quad c \,u^{\ast}(s) \geq v^{\ast}(s),
\;\;\; \forall s \in [0,A(\theta_{1})] \}
\end{eqnarray}
Then by the definition of  $c'$, we can find $ s_{1} \in [0, A(\theta_{1})],$ \; such that \; $c'\, u^{\ast}(s_{1}) =  v^{\ast}(s_{1})$.\\
We define now the following function
\begin{eqnarray*}h(s)=
\left\{\begin{array}{llll}
 c' \, u^{\ast}(s); &\text{ if }  s \in [0, s_{1}]\\
v^{\ast}(s)   ; &\text{ if } s \in [s_{1}, A(\theta_{1})]
\end{array}
\right.
\end{eqnarray*}
The properties of $u^{\ast}$ and $v^{\ast}$ imply that $h$ is monotonically decreasing and \\
$h(A(\theta_{1}))= 0$. Further, in virtue of (22) and (28), we
easily see that
\begin{eqnarray}\label{LLLL}
     -h'(s) \leq  \lambda (\beta \, \omega_{n-1})^{-2}(\sin \theta(s))^{2-2n}\int_{0}^{s} h(\xi)d\xi
   \end{eqnarray}
  for almost all  $s \in [ 0, A(\theta_{1})]$.\\
   Now, let $\mathfrak{h}$ be a radial function  defined in $B_{\theta_{1}(\lambda)}$ by
\begin{eqnarray}
     \mathfrak{h}(\theta) = h(A(\theta))
\end{eqnarray}
then, $\mathfrak{h}$ is an admissible function for the Rayleigh
quotient on $B_{\theta_{1}(\lambda)}$. From this we proceed exactly as in the
proof of  inequality (39), we obtain

\begin{eqnarray} \label{eq36}
\frac{\int_{B_{\theta_{1}(\lambda)}}|\nabla \mathfrak{h}|^{2}
dv_{g^{\star}}}{\int_{B_{\theta_{1}(\lambda)}} \mathfrak{h}^{2}  dv_{g^{\star}}}
\leq \lambda
\end{eqnarray}
and by the definition of $B_{\theta_{1}(\lambda)}$ it follows that
$\mathfrak{h}$ is an eigenfunction associated with $\lambda$ on
$B_{\theta_{1}(\lambda)}$. Then $\mathfrak{h}$ is a multiple of $v$ and so, from
the definition of $ h$ , it follows that $h = v^{\ast}$ and  $c'\,
u^{\ast}(s) = v^{\ast}(s)$ for $0 \leq s \leq s_{1}$ .  Since
$u^{\ast}(0) = v^{\ast}(0)$, then  $c' = 1$ and $ u^{\ast}(s) \geq
v^{\ast}(s)$ for all  $\; 0\leq s\leq  A(\theta_{1})$.
The proof of the lemma is thereby complete.\\\\
\textbf{Proof of theorem 1.1}
The normalization condition (1) is equivalent to
\begin{equation} \label{normalization1}
\int_0^{vol(D)} (u^{\ast}(s))^{p} d s = \int_0^{A(\theta_{1})}
(v^{\ast}(s))^{p} d s
\end{equation}
Since the functions $u^{\ast}$ and $v^{\ast}$ are nonnegative, and
$A(\theta_{1})\leq vol(D)$ (see \eqref{clm}), it is then clear that
\begin{equation} \label{normalization2}
\int_0^{A(\theta_{1})} (u^{\ast}(s))^{p} d s \leq
\int_0^{A(\theta_{1})} (v^{\ast}(s))^{p} d s
\end{equation}

We will first prove that $v^{\ast}(0)\ge u^{\ast}(0)$. Assume that $v^{\ast}(0) < u^{\ast}(0)$.In this case, $\exists \, \kappa  > 1$,
such that $\kappa\,v^{\ast}(0)= u^{\ast}(0)$. By Lemma3.1, we have
\begin{eqnarray}\label{A1A1A1A1}
\kappa\,v^{\ast}(s) \leq u^{\ast}(s) \quad \forall s \in [0,
A(\theta_{1})]
\end{eqnarray}
Therefore
\begin{equation*}
\kappa^p \int_{0}^{A(\theta_{1})}(v^{\ast}(s))^{p} d s \leq
\int_{0}^{A(\theta_{1})}(u^{\ast}(s))^{p} d s
\end{equation*}
Combining this inequality with \eqref{normalization2} leads to
$\kappa^p \le 1$, which is a contradiction.

Suppose now that $v^{\ast}(0) = u^{\ast}(0)$. From (43) and Lemma3.1, we obtain
\begin{eqnarray}
 \int_{0}^{vol(D)}(u^{\ast}(s))^{p} d s = \int_{0}^{A(\theta_{1})}(v^{\ast}(s))^{p} d s \leq \int_{0}^{A(\theta_{1})}(u^{\ast}(s))^{p} d s
\end{eqnarray}
This means $\displaystyle\int_{A(\theta_{1})}^{vol(D)}
(u^{\ast}(s))^{p} d s \leq 0$, and since $u^{\ast} > 0$ in $(0,
vol(D))$, we have $vol(D) = A(\theta_{1})$. Then $v = u^{\star}$,
and the statements of the theorem are evident.

Now, we treat the case $v^{\ast}(0) > u^{\ast}(0)$.

In this case $vol(D) > A(\theta_{1})$. Therefore,
$v^{\ast}(A(\theta_{1})) = 0$ and $u^{\ast}(A(\theta_{1})) > 0$.
Now, by the continuity of $v^{\ast}$ and $u^{\ast}$, we see that
$v^{\ast}(s) > u^{\ast}(s)$ in a neighborhood of 0, and there exists
$s_{1} \in (0,A (\theta_{1}))$ such that $v^{\ast}(s_{1}) =
u^{\ast}(s_{1})$. Choose $s_{1}$ to be the largest such number with
the additional property that $u^{\ast}(s)\leq v^{\ast}(s)$ \;for all
$s \in [0, s_{1}]$. By the definition of $s_{1}$, there is an
interval immediately to the right of $s_{1}$ on which $u^{\ast}(s) >
v^{\ast}(s)$. We will now show that $u^{\ast}(s) > v^{\ast}(s)$ for
all $s \in (s_{1}, A(\theta_{1})]$. If not, there exists $s_{2} \in
(s_{1}, A(\theta_{1}))$ such that $u^{\ast}(s_{2}) =
v^{\ast}(s_{2})$ and $u^{\ast}(s) > v^{\ast}(s)$ for all $s\in
(s_{1}, s_{2})$. In this case, we can define the function
\begin{eqnarray*}\varphi(s)=
\left\{\begin{array}{llll}
 v^{\ast}(s), &\text{ for } s \in [0, s_{1}] \cup [s_{2}, A(\theta_{1})],\\
u^{\ast}(s)    , &\text{ for } s \in [s_{1}, s_{2}]
\end{array}
\right.
\end{eqnarray*}
It follows from (22) and (28) that $\varphi$ satisfies
\begin{eqnarray}
     -\varphi'(s) \leq  \lambda \left(\beta \, \omega_{n-1}\right)^{-2}(\sin \theta(s))^{2-2n}\int_{0}^{s} \varphi(\xi)d\xi
   \end{eqnarray}
From $\varphi$ define a radial function in $B_{\theta_{1}(\lambda)}$ by
\begin{eqnarray}
     \Phi(\theta) = \varphi \big(A(\theta)\big)
   \end{eqnarray}
Then  $\Phi$ is an admissible function for the Rayleigh quotient on
$ B_{\theta_{1}(\lambda)}$. Using this fact we proceed exactly as in the proof
of  inequality (35), we obtain
\begin{eqnarray}
\frac{\int_{B_{\theta_{1}(\lambda)}}|\nabla \Phi|^{2}
dv_{g^{\star}}}{\int_{B_{\theta_{1}(\lambda)}} \Phi^{2} dv_{g^{\star}}} \leq
\lambda
\end{eqnarray}

   It follows that the Rayleigh quotient  of $\Phi$ is
equal to $\lambda$ therefore $\Phi$ is an eigenfunction for
$\lambda$, Consequently, $u^{\ast}= v^{\ast}$ and so $u^{\ast}(s) =
v^{\ast}(s)$ in $[s_{1}, s_{2}]$ which contradicts the maximality of
$s_{1}$ and hence completes the proof of the theorem.\\

\textbf{Proof of Corollary1.2: Chiti's Reverse H\"older Inequality}\\
For $p>0$, we choose $v$   so that the condition (1) is satisfied.
Now, if we extend the function $v^{\ast}$  by  zero in
$[A(\theta_{1}), vol(D) ]$, we obtain
\begin{eqnarray} \label{eq3.50}
\int_0^{s} (u^{\ast}(\eta))^{p} d\eta \leq \int_0^{s}
(v^{\ast}(\eta))^{p} d \eta \;\;\;\;\;\;\forall s \in [0, vol(D)]
\end{eqnarray}
To see \eqref{eq3.50}, we note that Theorem1.1  implies the
following:

If $s \in [0, s_{1}]$, then
\[u^{\ast}(\eta)  \leq v^{\ast}(\eta), \;\;\; \forall \eta \in [0, s]\]
hence
$$\int_{0}^{s}(u^{\ast}(\eta))^{p} d\eta \leq \int_{0}^{s}(v^{\ast}(\eta))^{p}d\eta$$

If $s \in [s_{1}, vol(D)]$, then
\begin{eqnarray*} \int_{0}^{s}(u^{\ast}(\eta))^{p} d\eta
&=& \int_{0}^{vol(D)}(u^{\ast}(\eta))^{p} d\eta - \int_{s}^{vol(D)}(u^{\ast}(\eta))^{p} d\eta \notag \\
&\leq& \int_{0}^{vol(D)}(v^{\ast}(\eta))^{p} d\eta - \int_{s}^{vol(D)}(v^{\ast}(\eta))^{p} d\eta  \notag \\
&=& \int_{0}^{s}(v^{\ast}(\eta))^{p} d\eta
\end{eqnarray*}
We complete the argument using the following result
\begin{lemma} \cite{HLP}
 Let  $R, p, q$  be real numbers such that $0<p\leq q$, \; $R>0$; and  $f,$ $g$  real
functions in $L^{q}([0,R])$. If the decreasing rearrangements of $f$
and $g$ satisfy the following inequality:
$$\int_{0}^{s} \left(f^{\ast}\right)^{p}dt\leq \int_{0}^{s}\left(g^{\ast}\right)^{p}dt, \quad \forall \; s\in [0,R], $$  \\
then
$$\int_{0}^{R}f^{q}dt\leq \int_{0}^{R}g^{q}dt.$$
\end{lemma}
From this, it is clear, that for all $q \geq p$,
\begin{eqnarray}
\int_{0}^{vol(D)}(u^{\ast}(\eta))^{q} d\eta \leq
\int_{0}^{vol(D)}(v^{\ast}(\eta))^{q} d\eta =
\int_{0}^{A(\theta_{1})}(v^{\ast}(\eta))^{q} d\eta
\end{eqnarray}
Finally, combining this inequality and equality (1), we obtain the desired result.\\
Now assume that we have equality in (4), from the normalization for
the function $v$ given in (1), we deduce that for all $p>0,$
$\int_{D} u^{p} dv_{g} =\beta \, \int_{B_{\theta_{1}(\lambda)}} v^{p}
dv_{g^{\star}}.$ Hence $vol(D) = \beta vol(B_{\theta_{1}(\lambda)}),$ and since
$D^{\star}$ and $B_{\theta_{1}(\lambda)}$ are geodesic balls of $\mathbb{S}^{n}$
centered both at the north pole with the same volume, it yields that
$D^{\star}= B_{\theta_{1}(\lambda)}.$ By hypothesis $\lambda$ is the fundamental
eigenvalue
 of $B_{\theta_{1}(\lambda)},$ hence of $D^{\star},$ thus we obtain that $\lambda_{1}(D)= \lambda_{1}(D^{\star})= \lambda$
 and this is possible if and only if the triplet $(M,D, g)$ is isometric to
the triplet $(\mathbb{S}^{n}, D^{\star},g^{\star}),$  ( one can see
Theorem 5 in \cite{BM}). The proof of Corollary1.2 is thereby
complete.
\subsection{Saint-Venant Theorem}
Let $(M, g)$ be a compact  Riemannian manifold of dimension $n$,
without boundary satisfying $R(M, g)\geq  n-1$, and let $D$ a
 bounded connected domain of $M$ with smooth boundary. We are interested in the
following geometric quantity
 \begin{eqnarray}
T(D) = \int_{D} w dv_{g}
\end{eqnarray}
 where $w$ is is the  smooth solution of the boundary value problem of Dirichlet-Poisson type
  \begin{eqnarray} \label{ss}
\Delta_{g} w + 1 & = & 0 \mbox{ in } D, \nonumber \\
w & = & 0 \mbox{ on }\partial D
\end{eqnarray}
 The geometric quantity $T(D)$ is
called the ``torsional rigidity of $D$'', and it is customary to
call the solution $w$ of \eqref{ss} the warping function.
In theorem1.3 we will give  explicit upper bounds for the
torsional rigidity $T(D)$ which amounts to a version of the
Saint-Venant Theorem for compact manifolds. Let $C_{0}^{\infty}(D)$
denotes the space of $C^{\infty}$ functions with compact support in
$D$. We define the Sobolev space $H_{0}^{1}(D)$ as the closure of
$C_{0}^{\infty}(D)$ in $H^{1}(M)$ the space of square integrable
functions with  a square integrable weak derivatives. The
variational formulation for $T(D)$ is given by
 \begin{eqnarray}
\frac{1}{T(D)} = \inf \left\{ \Phi(f) = \frac{\int_{D}|\nabla u|^{2}
dv_{g}}{(\int_{D} u dv_{g})^{2}},\;\;\; f \in H_{0}^{1}(D), \;\;\;
f\ncong 0  \right\}
\end{eqnarray}
Indeed, by the scaling property of the functional, $\Phi(c f) =
\Phi(f)$ for all $c > 0$, one can reformulate the minimizing problem
of the functional $\Phi$ as a minimizing problem of the functional
$\int_{D}|\nabla f|^{2} dv_{g}$ subject to the constraint  $\int_{D}
f dv_{g}= 1$.  By the above mentioned Lagrange multipliers theorem,
this gives the existence of the Lagrange multiplier $\rho$, such
that for any $h \in H_{0}^{1}(D)$
\begin{eqnarray}
\int_{D}\langle \nabla f, \nabla h  \rangle dv_{g}= \rho \int_{D} h
dv_{g}
\end{eqnarray}
Hence, $f$ is a weak solution of the equation
\begin{eqnarray}
\Delta f & = -\rho& 0 \mbox{ in } D, \nonumber \\
f& = & 0 \mbox{ on }\partial D
\end{eqnarray}
By standard regularity results, $f$ is unique and smooth. Then $w =
\rho^{-1}f$ is also a critical point of $\Phi$. Finally the fact
that $\Phi (w) = \frac{1}{T(D)}$ proves the equality.\\

 We will now give the proof  of the Saint-Venant
Theorem.\\

\textbf{Proof of Theorem1.3}\\
 The proof follows the steps
of Talenti's method (one can see \cite{Talenti}) tailored to our
setting. For $0 <t\leq m = \sup \{w(x),\;\; x\in D\} $, we define
the set
\begin{eqnarray}
D_{t} = \{x \in D;\;\;\;\; w(x)>t\}
\end{eqnarray}
and  introduce the following functions
\begin{eqnarray}
V(t) = \int_{D(t)}dv_{g}, \;\;\;\;\;\;\;\;\;\Psi(t) =
\int_{D(t)}|\nabla w|^{2}dv_{g}
\end{eqnarray}
  The smooth co-area formula gives
\begin{eqnarray}
   V(t) =  \int_{t}^{m}\int_{\partial D_{\tau}}\frac{1}{|\nabla w|}  d\sigma d\tau,\;\;\;\;\;\;  \Psi(t) =  \int_{t}^{m}\int_{\partial D_{\tau}}|
   \nabla w|  d\sigma d\tau
\end{eqnarray}
Then differentiating (58) with respect to $t,$ one obtain
\begin{eqnarray}
   V^{\prime}(t) =  -\int_{\partial D_{t}}\frac{1}{|\nabla w|}d\sigma ,\;\;\;\;\;\;  \Psi^{\prime}(t) =  -\int_{\partial D_{t}}|\nabla w|  d\sigma
\end{eqnarray}

Let $w^{\ast}$ be the inverse function of $V,$  $w^{\star}$  the
radial function defined in $D^{\star}$ by $w^{\star}(\theta) =
w^{\ast}(A(\theta)),$ and
 \begin{eqnarray}
V_{\star}(t) = \int_{D_{t}^{\star}}dv_{g^{\star}},
\;\;\;\;\;\;\;\;\;\;\Psi_{\star}(t) = \int_{D_{t}^{\star}}|\nabla
w^{\star}|^{2}dv_{g^{\star}}
\end{eqnarray}

Since  $w^{\star}$ is a radial  decreasing function, its
  level sets $D^{\star}_{t}$ are  geodesic balls with
radius $r(t)  =w^{\star -1}(t)$. Therefore
\begin{eqnarray}
V_{\star}(t) = \omega _{n-1}\int_{0}^{r(t)}\sin^{n-1}\theta d\theta
\end{eqnarray}
and
\begin{eqnarray}
 \Psi_{\star}(t)= \omega _{n-1}\int_{0}^{r(t)}(\frac{d w^{\star}}{d \theta})^{2}\sin^{n-1}\theta d\theta
\end{eqnarray}
Differentiating with respect to $t$, we get
\begin{eqnarray}
V'_{\star}(t) = -|\nabla
w^{\star}|^{-1}(t)\omega_{n-1}\sin^{n-1}r(t) = - |\nabla
w^{\star}|^{-1}(t)\int_{\partial D^{\star}_{t}} d\sigma
\end{eqnarray}
and
\begin{eqnarray}
 \Psi'_{\star}(t) = -|\nabla w^{\star}|(t)\omega _{n-1}\sin^{n-1}r(t) = - |\nabla w^{\star}|(t)\int_{\partial D^{\star}_{t}} d\sigma
\end{eqnarray}
  On one hand multiplying  (64) by (65), we obtain
\begin{eqnarray}
\left(\int_{\partial D^{\star}_{t}} d\sigma\right)^{2}
=V'_{\star}(t) \Psi'_{\star}(t)
\end{eqnarray}
In the other hand the Cauchy-Schwartz inequality gives
\begin{eqnarray}
\left(\int_{\partial D_{t}} d\sigma\right)^{2}  \leq V'(t) \Psi'(t)
\end{eqnarray}
   From the fact that $V_{\star}(t) = \beta^{-1}\, A(r(t))= \beta^{-1}\, V(t)$, we have
    $V'(t) = \beta \, V'_{\star}(t)$, then we use  Gromov's isoperimetric inequality (14), (66) and (67) to obtain
\begin{eqnarray}
\beta \, \Psi'_{\star}(t) \geq \Psi'(t)
 \end{eqnarray}
 Now by integrating this inequality from $0$ to $m$ and using the fact that $\beta \, \Psi_{\star}(m)= 0 =  \Psi(m),$
 we get
 \begin{eqnarray}
\beta \, \int_{D^{\star}}|\nabla w^{\star}|^{2}dv_{g^{\star}} =\beta
\, \Psi_{\star}(0) \leq \Psi(0) = \int_{D}|\nabla w|^{2}dv_{g}
\end{eqnarray}
 Next, the fact that
\begin{eqnarray}
\beta \, \int_{D^{\star}} w^{\star}dv_{g^{\star}}  = \int_{D} w
dv_{g}
\end{eqnarray}
and the  variational formulation given in (54), it yields
\begin{eqnarray}
\frac{1}{T(D)} = \frac{\int_{D}|\nabla w|^{2}dv_{g}}{(\int_{D} w
dv_{g})^{2}} \geq \frac{1}{\beta} \, \frac{ \int_{D^{\star}}|\nabla
w^{\star}|^{2}dv_{g^{\star}}}{(\int_{D^{\star}}
w^{\star}dv_{g^{\star}})^{2}}\geq \frac{1}{\beta \, T(D^{\star})}
\end{eqnarray}
In the case of  equality, going back to (68), which we integrate,
then, we use the facts that
 $\beta \, \Psi_{\star}(m)= \Psi(m)$
and $\beta \, \Psi_{\star}(0) = \Psi(0),$  we prove that $\beta\,
\Psi'_{\star}(t)= \Psi'(t)$ for all $0 <t\leq m .$ Next, we use the
last equality in (66) and (67), we get
\begin{eqnarray}
vol_{n-1}(\partial D^{\star}_{t}) = \int_{\partial D^{\star}_{t}}
d\sigma = \int_{\partial D_{t}} d\sigma = vol_{n-1}(\partial D_{t})
\end{eqnarray}
to achieve the proof, we apply once again Gromov's isoperimetric
inequality.\\

In the sequel, we will  give a comparison theorem for the warping
function in the case of a smooth compact Riemannian manifold. This
theorem is based on a method of Talenti (\cite{Talenti}).

\begin{theorem}

Let $(M, g)$ be a compact  Riemannian manifold of dimension $n$,
without boundary satisfying $R(M, g)\geq  n-1$, and  $D$ a
 bounded connected domain of $M$ with smooth boundary.
Let $v$ be the radial function defined in $D^{\star}$ by
\begin{eqnarray}
v(\theta) = \int_{\theta}^{\theta_{0}}
\left(\int_{0}^{\delta}\sin^{n-1}\tau  d\tau \right)\sin^{1-n}
\delta d\delta
\end{eqnarray}
 which is the solution to the problem
 \begin{eqnarray}
\Delta_{g} v + 1 & = & 0 \mbox{ in } D^{\star}, \nonumber \\
v & = & 0 \mbox{ on }\partial D^{\star}
\end{eqnarray}
Then $w^{\star}$, the symmetrized function  of $w$ satisfies
 \begin{eqnarray}
w^{\star} \leq v\;\;\;\;\;\;\;\;\;\;a.e \;\;in\;\; D^{\star}
\end{eqnarray}
 We obtain the equality if and only if the triplet $(M,D,g)$ is isometric to the triplet $(\mathbb{S}^{n},D^{\star},g^{\star})$
\end{theorem}
\textbf{Proof}\\

Let $f$ be a test function in the weak formulation of our problem
defined on $D$  by
\begin{eqnarray} \label{test} \; f(x) =
\left\{\begin{array}{llll}
w(x)- t, &\text{ if } w(x)>t \\
0  &\text{otherwise }
\end{array}
\right.
\end{eqnarray}
where $ 0\leq t < m$.  We introduce the function defined by
\begin{eqnarray}
\Psi(t) = \int_{D_{t}}|\nabla  w|^{2}dV_{g}
\end{eqnarray}
Then
\begin{eqnarray}
     \Psi(t)  = \int_{w> t}\, (w- t)\,  dV_{g}
\end{eqnarray}
The function $\Psi$ is decreasing  in $t$ then, for $h > 0$ we have
\begin{eqnarray*}
\frac{\Psi(t) - \Psi(t + h) }{h} = \int_{w> t + h} \, dV_{g} +
\int_{t < w\leq t+h}\, \left(\frac{w- t}{h}\right)  \, dV_{g}
\end{eqnarray*}
Letting $h$ go to zero, we obtain, for the right derivative of $\Psi
(t)$,
\begin{eqnarray}
-\Psi'_{+}(t) = \int_{w > t} dV_{g} \quad \text{a.e.} \quad t > 0
\end{eqnarray}
The same computation gives the same equality for the left derivative
of $\Psi(t)$. Therefore
\begin{eqnarray}\label{ineq1001}
0 \leq -\Psi'(t)= V(t)
\end{eqnarray}
Next, we use the Cauchy-Schwarz inequality
\begin{eqnarray*}
 \left(\frac{1}{h} \int_{t < w \leq t+h} |\nabla w| dV_{g}\right)^{2}
\leq \left(\frac{1}{h} \int_{t < w \leq t+h} |\nabla w|^{2}
dV_{g}\right)\left(\frac{1}{h} \int_{t < w \leq t+h}dV_{g} \right)
\end{eqnarray*}
Thus, letting $h\rightarrow 0$ and using \eqref{ineq1001}, we obtain
\begin{eqnarray}\label{eq89}
\left(-\frac{d}{d t}\int_{f> t}|\nabla w |dV_{g} \right)^{2} \leq
(-V'(t)V(t))
\end{eqnarray}
Using the co-area formula, it yields
 \begin{eqnarray} \label{eq88}
-\frac{d}{dt}(\int_{w> t}|\nabla w|dV_{g}) = \int_{\partial
D_{t}}d\sigma , \quad \text{a.e.} \quad t> 0
\end{eqnarray}
Next by M.Gromov's isoperimetric inequality (14), we have
\begin{eqnarray}
\int_{\partial D_{t}}d\sigma \geq \beta \, \int_{\partial
D^{\star}_{t}} d\sigma =\beta \, \omega _{n-1}\sin^{n-1}r(t) =
A'(r(t))
\end{eqnarray}
Combining this with \eqref{eq88} and \eqref{eq89}, we get
\begin{eqnarray}\label{eq1005}
-V'(t) V(t)   \geq   \left(A'(r(t))\right)^{2}
\end{eqnarray}
 for almost every $t$ in $(0, m)$. Using the fact that $A(r(t)) = V(t)$ and $r^{-1}(\theta) = w^{\star}(\theta),$ we get
\begin{eqnarray} \label{last-ineq}
-w^{\star '}(\theta) \leq \frac{A(\theta)}{A'(\theta)}
\end{eqnarray}
Now, for $\theta \in (0,\theta_{0})$, integrating this inequality
from $\theta$  to $\theta_{0},$ we obtain
\begin{eqnarray}
w^{\star }(\theta) \leq
\int_{\theta}^{\theta_{0}}\frac{A(\tau)}{A'(\tau)} d\tau =
v(\theta)
\end{eqnarray}
which is the desired result.\\
 Now assume that we have equality in (7),
integrating this equality, we get
\begin{eqnarray}
T(D) = \int_{D} w dV_{g}=\beta \, \int_{D^{\star}} w^{\star}
dV_{g^{\star}} =\beta \, \int_{D^{\star}} v dV_{g^{\star}}
\end{eqnarray}
Finally by applying the Saint-Venant Theorem, we deduce that the
triplet $(M,D,g)$ is isometric to the triplet
$(\mathbb{S}^{n},D^{\star}, g^{\star}),$ which completes the proof
of the theorem.\\
\begin{remark} As a consequence of the result above, we obtain after integrating
inequality (75):\\
\begin{eqnarray}
T(D) =\int_{D} w dv_{g}= \beta \int_{D^{\star}} w^{\star}
dv_{g^{\star}}\leq \beta\int_{D^{\star}} v dv_{g^{\star}}=
T(D^{\star})
\end{eqnarray}
\end{remark}
Which is the result of Theorem1.3.

E-mails: ngamara@taibahu.edu.sa\\
         \hspace*{15mm}\quad\quad  hasnaouiabdelim9@gmail.com\\
         \hspace*{15mm} \quad \quad akremmakni@gmail.com\\

Addresses: College of Sciences, Taibah University, Kingdom of Saudi Arabia.\\
           \hspace*{17mm}\quad \quad \quad  University Tunis El Manar, FST, Mathematics Department, Tunisia.\\


\section*{References}
\begin{thebibliography}{xxx}
\bibitem{krahn} Krahn, E.: \H{U}ber eine von Rayleigh formulierte Minmaleigenschaft des Kreises, Math.\ Ann.\ \textbf{94}, 97--100  (1925).

\bibitem{faber} Faber, C.: Beweiss, dass unter allen homogenen Membrane von gleicher Fl\H{a}che und gleicher Spannung die kreisf\H{o}rmige die tiefsten Grundton gibt. Sitzungsber.-Bayer. Akad. Wiss., Math.-Phys. Munich., pp. 169--172 (1923).
\bibitem{PR1} Payne, L.~E., Rayner, M.~E.: Some isoperimetric norm bounds for solutions of the Helmholtz equation, Z.\ Angew.\ Math.\ Phys.\ \textbf{24}, 105--110  (1973).
\bibitem{PR2} Payne, L.~E., Rayner, M.~E.: An isoperimetric inequality for the first eigenfunction in the fixed membrane problem, Z.\ Angew.\ Math.\ Phys.\ \textbf{23}, 13--15  (1972).
\bibitem{KJ1} Kohler-Jobin, M.T.: Sur la premi\`ere fonction propre d'une membrane: une extension \`a $N$ dimensions de l'in\'egalit\'e isop\'erim\'etrique de Payne-Rayner, Z.\ Angew.\ Math.\ Phys.\ \textbf{28},  1137--1140  (1977).
\bibitem{Ch1} Chiti, G.: A reverse H\"older inequality for the eigenfunctions of linear second order elliptic operators, Z. Angew. Math. Phys. \textbf{33}, 143-148  (1982).
\bibitem{H.Has}  Hasnaoui, A.: On the Problem of Queen Dido for a Wedge like membrane and a Compact Riemannian manifold with lower Ricci curvature bound, PHD January 2014, University Tunis El Manar, Tunisia.

\bibitem{Cheng} S.Y. Cheng: Eigenvalue comparison theorems and its geometric aplications, Math. Z. 143 pp 289-297,(1975).
\bibitem{BM} B\'erard, P.,  Meyer, D.: In\'egalit\'es isop\'erim\'etriques et applications, Ann. Scient. Ec. Norm. Sup. ´
(4), 15 (1982) 513–542.
\bibitem{Gam} N. Abdelmoula Gamara: Sym\'etrisation d'in\'equations elliptiques et applications g\'eom\'etriques, Math. Z. 199 pp 181-190,(1988).
\bibitem{LingLu} Ling, J., Lu, Z.: Bounds of eigenvalues on Riemannian manifolds. Trends in partial differential equations, 241–264, Adv. Lect. Math. (ALM), 10, Int. Press, Somerville, MA (2010).
\bibitem{Saint-Venant} De Saint-Venant, B.: M\'emoire sur la torsion des prismes. M\'emoir. pres. divers. savants, Acad. Sci. \textbf{14}, 233--560 (1856).
\bibitem{polya48} P\'olya, G.: Torsional rigidity, principal frequency, electrostatic capacity, and symmetrization, Quart.\ of Appl.\ Maths \textbf{6}, 267--277 (1948).
\bibitem{Makai} Makai, E.: A proof of Saint-Venant's theorem on torsional rigidity, Acta Math.\ Acad.\ Sci.\ Hungar. \textbf{17}, 419--422  (1966).
\bibitem{Bandle} Bandle, C.: Isoperimetric inequalities and applications, Monographs and Studies in Mathematics, 7, Pitman (Advanced Publishing Program), Boston, Mass.-London (1980).
\bibitem{PolyaSzego} P\'olya, G., Szeg\H{o}, G.: Isoperimetric Inequalities in Mathematical Physics. Princeton University Press (1951).
\bibitem{RatzkinCaroll2011} Carroll, T., Ratzkin, J.: Interpolating between torsional rigidity and principal frequency, J.\ Math.\ Anal.\ Appl.\ \textbf{379}, 818--826  (2011).
\bibitem{Iversen} Iversen, M.: Torsional rigidity of a radially perturbed ball, Oberwolfach Report \textbf{33}, 31--33 (2012).
\bibitem{vdb0} Van den Berg, M.:  Estimates for the torsion function and Sobolev constants, Potential Anal.\ 36 (2012) 607--616.
\bibitem{RatzkinCaroll2014} Carroll, T., Ratzkin, J.: A reverse H¨older inequality for extremal Sobolev functions, Potential Anal, DOI 10.1007/s11118-014-9433-6.
\bibitem{Payne1}  Payne, L.~E.: Some comments on the past fifty years of isoperimetric inequalities, Inequalities (Birmingham, 1987), 143--161, Lecture Notes in Pure and Appl. Math., 129, Dekker, New York, 1991.
\bibitem{Payne0} Payne, L.~E.:  Some isoperimetric inequalities in the torsion problem for multiply connected regions, Studies in mathematical analysis and related topics: Essays in honor of Georgia P\'olya, Stanford Univ. Press, Stanford, Calif.,  pp. 270--280, 1962.
\bibitem{PayneWeinbergerTorsion} Payne, L.~E., Weinberger, H.~F.: Some isoperimetric inequalities for membrane frequencies and torsional rigidity, J.\ Math.\ Anal.\ Appl.\ \textbf{2}, 210-216  (1961).
\bibitem{AshbaughBenguriaSphere} Ashbaugh, M. S., Benguria, R. D.: A sharp bound for the ratio of the first two Dirichlet eigenvalues of a domain in a hemisphere of Sn. Trans. Amer. Math. Soc. 353, 1055-–1087 (2001).
\bibitem{ChavelFeldman} Chavel, I., Feldman, E. A.: Isoperimetric inequalities on curved surfaces. Adv. in Math. 37, 83–-98 (1980).
\bibitem{Gromov} Gromov, M.: Paul Levy's isoperimetric inequality. Appendix C in Metric structures for Riemannian and non-Riemannian spaces. Based on the 1981 French original. With appendices by M. Katz, P. Pansu and S. Semmes. Translated from the French by Sean Michael Bates. Progress in Mathematics, 152. Birkh\H{a}user Boston, Inc., Boston, MA, 1999.
\bibitem{Benguria2011} Benguria, R. D.:
Isoperimetric inequalities for eigenvalues of the Laplacian. Entropy
and the quantum II, 21--60, Contemp. Math., \textbf{552}, Amer.
Math. Soc., Providence, RI (2011).
\bibitem{Ch2} Chiti, G.: An isoperimetric inequality for the eigenfunctions of linear second order elliptic operators, Boll. Un. Mat. Ital. A (6) \textbf{1}, 145-151  (1982).
\bibitem{HLP} Hardy, G.~H., Littlewood, J.~E., P\'{o}lya, G.: Some simple inequalities satisfied by convex functions, Messenger Math. 58,
145--152 (1929). Reprinted in Collected Papers of G. H. Hardy, Vol.II, London Math. Soc., Clarendon Press: Oxford, pp. 500--508 (1967).

\bibitem{Talenti} Talenti, G.: Elliptic equations and rearrangements, Ann.\ Scuola Norm.\ Sup.\ Pisa Cl.\ Sci.\ \textbf{3}, 697--718  (1976).



%%%%%%%%%%%%%%%%%%%%%%%%%%%%%%%%%%%%%%%%%%%%%%%%%%%%%%%%%%5




\end{thebibliography}
\end{document}